\documentclass[conference]{IEEEtran}

\usepackage[utf8]{inputenc}
\usepackage{hyperref}
\usepackage{graphicx}
\usepackage{amsmath} 
\usepackage{centernot}
\usepackage[dvipsnames]{xcolor}
\DeclareMathOperator{\diag}{diag}
\DeclareMathOperator{\rank}{rank}
\newcommand*{\defeq}{\stackrel{\text{def}}{=}}

\title{Pitfalls of Zero Voltage Values in \\Optimal Power Flow Problems}

\author{\IEEEauthorblockN{Frederik Geth}
\IEEEauthorblockA{
\textit{GridQube}\\
Brisbane, Australia \\
frederik.geth@gridqube.com}
}

\begin{document}

\maketitle

\newcommand{\ivr}{IVR}
\newcommand{\svrone}{SVR-1}
\newcommand{\svrtwo}{SVR-2}
\newcommand{\swrone}{SWR-1}
\newcommand{\swrtwo}{SWR-2}

\begin{abstract}
The existence of strictly positive lower bounds on voltage magnitude is taken for granted in optimal power flow problems. 
Nevertheless, it is not possible to rely on such bounds for a variety of real-world network optimization problems.
This paper discusses a few issues related to 0\,V assumptions made during the process of deriving optimization formulations in the current-voltage, power-voltage and power-lifted-voltage variable spaces. 
The differences between the assumptions are illustrated for a 2-bus 2-wire test case, where the feasible sets are visualized. 
A nonzero relaxation gap is observed for the canonical multiconductor nonlinear power-voltage formulation. 
A zero gap can be obtained for the branch flow model semi-definite relaxation, using newly proposed valid equalities.
\end{abstract}


\section{Introduction}

In the context of ac transmission system optimization, it has generally been established \cite{kardos2018complete} that the nonlinear nonconvex power-voltage formulation of the circuit physics (in \textbf{p}olar coordinates, `SV\textbf{P}') is the most scalable, when used with nonlinear programming solvers such as \textsc{Ipopt} \cite{ipopt}, MIPS  \cite{5491276} or \textsc{Beltistos} \cite{KARDOS2022108613}.
Other exact formulations, e.g. the \textbf{r}ectangular current-voltage  (IV\textbf{R}) and the rectangular power-voltage formulation (SVR) result in overall slower computations \cite{kardos2018complete}. 

The canonical optimal power flow (OPF) problem specification furthermore assumes the existence of a number of bounds, and postulates a generation cost minimization objective. 
In general authors assume that voltage magnitude has a strictly positive lower bound, thereby excluding the 0\,V  from the feasible set.
Nevertheless, many real-world network optimization problems may not have this property of strictly positive lower bounds, as 1) the values may be missing entirely; or 2)  0\,V solutions can not be excluded a-priori.
In the author's experience, such contexts are  not far-fetched:
\begin{itemize}
    \item Short circuit situations, e.g. in the context of (security-constrained) OPF with explicit representation of short circuit conditions \cite{vanacker2015,en14082160}.
    \item Polynomial chaos  OPF \cite{7587824,vanacker2022} has no magnitude lower bound for all $k>1$ (i.e. there is a valid lower bound on the mean, but the 0\,V solution cannot be excluded from the variance).
    \item Harmonic OPF \cite{geth2022} has no magnitude lower bound for the true harmonics $h>1$ (e.g. voltage magnitude at $h=3$, i.e., the triplen harmonic, can perfectly be 0).
    \item Four-wire OPF \cite{ARAUJO2013632,claeys2022} has no voltage magnitude lower bound for the neutral terminal\footnote{Even in multi-earthed networks, a behind-the-meter network may have neutral voltage rise, making Kron's reduction of the neutral  problematic.}. 
    \item OPF with HVDC lines has no voltage lower bound on the nodes for a metallic return~\cite{chandra2022}.
\end{itemize}
In grounded three-phase four-wire networks, the neutral current may be (close to) zero due to balanced loading, so excluding the 0\,V neutral solution may be impossible.
Allowing 0\,V values causes real algorithmic problems for SVR \cite{claeys2022}:
\begin{itemize}
    \item  0\,V solutions were naturally found by \textsc{Ipopt};
    \item  the solution obtained was very sensitive to variable initialization, i.e. different nodes ended up at 0\,V;
    \item it led to convergence problems, requiring much larger computation time than the \ivr{} form;
    \item solutions were likely to be \emph{strictly locally optimal}.
\end{itemize}
Similarly, polynomial chaos OPF paper \cite{vanacker2022} observed a significant computational speedup using \ivr{} over SVR. 

Therefore, this paper showcases a number of pitfalls related to 0\,V assumptions and solutions in the feasible sets of OPF in different variable spaces. 
The  0\,V solutions lead to problems in the lifting of Kirchhoff's current law (KCL) to power variables. 
Furthermore, loads and generators are often assumed to be connected between a phase and the neutral at 0\,V, and dropping this assumption requires attention to detail.

The paper is structured as follows.
First we derive a variety of formulations in a consistent notation, which are compared numerically next on a 2-bus 2-wire (explicit neutral) power flow case study\footnote{Extensive numerical results with 0\,V issues available in  \cite{vanacker2022,claeys2022}.}. 
Finally we present some conclusions and directions of future work.

\section{Derivation of Formulations}

\begin{figure}[tbh]
  \centering
    \includegraphics[width=0.75\columnwidth]{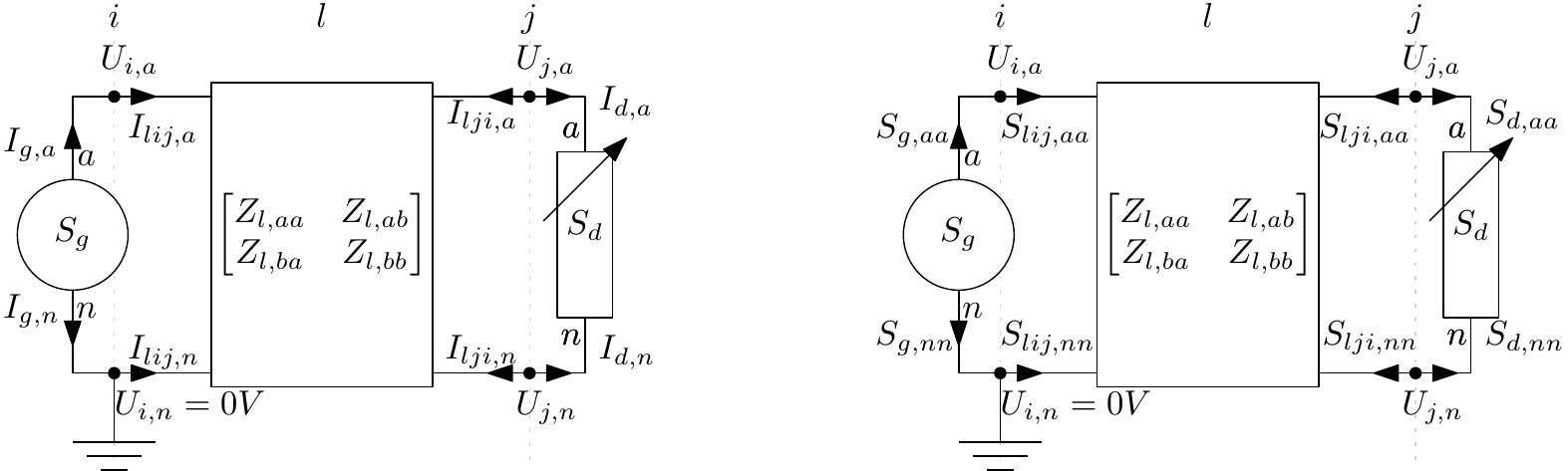}
  \caption{Single-grounded two-wire two-bus system with current (left) and power (right) as flow variables.}  
  \label{fig_SU-IV}
\end{figure}

Fig.~\ref{fig_SU-IV} illustrates the definitions of the current, voltage and power variables and indices used in this paper. 
We will explore the  following three-phase \emph{four-wire} a.k.a. \emph{unbalanced} rectangular formulations:
\begin{itemize}
    \item NLP unbalanced IVR  form, presented in \cite{claeys2022};
    \item NLP SVR form, with extensions to avoid violation of KCL in loads/generators, as discussed in \cite{claeys2022};
    \item Branch flow model (BFM) SDP relaxation \cite{gan2014},   extended to avoid violation of KCL in load/generators.  
\end{itemize}
We derive the n-conductor expressions in a consistent notation in vector-matrix variables, but \emph{illustrate} for a 2-wire system with  phase $a$ and an explicit neutral $n$. 
We ignore branch shunt terms for brevity, and refer to \cite{geth2021a} for a treatment of 4-wire OPF with asymmetric $\Pi$-model branches, and bus shunts.
 
\subsection{Current-Voltage (\ivr)}
The complex voltage vector $\mathbf{U}_i$ at bus $i$, a.k.a. phasor, stacks the values of the nodes at buses $i$,
\begin{IEEEeqnarray}{RL?s}
\mathbf{U}_i 
\defeq 
\begin{bmatrix}
U_{i,a} &
U_{i,n}
\end{bmatrix}^{\text{T}}.
\end{IEEEeqnarray}
For a two-wire set of electromagnetically coupled conductors, we end up with a $2\times 2$ impedance matrix, capturing the self-impedance (diagonal) and mutual impedance (off-diagonal). By taking the \emph{matrix inverse}, we obtain the admittance,
\begin{IEEEeqnarray}{RL?s}
\mathbf{Z}_l 
\defeq
\mathbf{R}_l + j \mathbf{X}_l 
\defeq
\begin{bmatrix}
Z_{l,aa} & Z_{l,an} \\
Z_{l,na} & Z_{l,nn} 
\end{bmatrix}, \mathbf{Y}_l \defeq (\mathbf{Z}_l)^{-1}.
\end{IEEEeqnarray}
Note that, based on the physics (Carson/Pollaczek eq.), we expect that $\mathbf{R}_l$ and $\mathbf{X}_l$  are symmetric and positive definite. 
We define the complex current vector in branch $l$ in the direction of $i$ to $j$ to stack the individual values for each of the wires, 
\begin{IEEEeqnarray}{RL?s}
\mathbf{I}_{lij} \defeq 
\begin{bmatrix}
I_{lij,a}&
I_{lij,n}
\end{bmatrix}^{\text{T}}, \quad
\mathbf{I}_{lji} \defeq
\begin{bmatrix}
I_{lji,a}&
I_{lji,n}
\end{bmatrix}^{\text{T}}.
\end{IEEEeqnarray}
The branch currents satisfy,
\begin{IEEEeqnarray}{RL?s}
\mathbf{I}_{lij} + \mathbf{I}_{lji}= 0. \label{eq_kcl_branch}
\end{IEEEeqnarray}
The multiconductor ac Ohm's law is,
\begin{IEEEeqnarray}{RL?s}
\mathbf{U}_j = \mathbf{U}_i - \mathbf{Z}_l  \mathbf{I}_{lij} . \label{eq_ohms_impedance}
\end{IEEEeqnarray}
The power flow into the branch is,
\begin{IEEEeqnarray}{RL?s}
\mathbf{S}_{lij} 
\defeq 
\begin{bmatrix}
S_{lij,aa}\\
S_{lij,nn}
\end{bmatrix}
\defeq
\mathbf{U}_i \circ (\mathbf{I}_{lij})^*
=
\begin{bmatrix}
U_{i,a}I_{lij,a}^*\\
U_{i,n}I_{lij,n}^*
\end{bmatrix} \label{eq_line_power_def},
\end{IEEEeqnarray}
where `$\circ$' indicates the element-wise product.
The current vector for a generator is,
\begin{IEEEeqnarray}{RL?s}
\mathbf{I}_{g} 
\defeq 
\begin{bmatrix}
I_{g,a}&
I_{g,n}
\end{bmatrix}^{\text{T}},
\end{IEEEeqnarray}
and by convention flows from the generator to the bus it is connected to, in each of its wires.
Note that  conservation of current applies through the generator, 
\begin{IEEEeqnarray}{RL?s}
 \mathbf{1}^\text{T} \mathbf{I}_{g}  = 0. \label{eq_kcl_gen}
\end{IEEEeqnarray}
The power  from the generator through each of the wires is, 
\begin{IEEEeqnarray}{RL?s}
\mathbf{S}_{g} 
\defeq 
\begin{bmatrix}
S_{g,aa} \\
S_{g,nn}
\end{bmatrix}
\defeq
\mathbf{U}_{i}
\circ
\mathbf{I}_{g}^*
=
\begin{bmatrix}
U_{i,a} I_{g,a}^*\\
U_{i,n} I_{g,n}^*
\end{bmatrix} \label{eq_gen_power_def}.
\end{IEEEeqnarray}
For loads, the current flows from the bus to the load, 
\begin{IEEEeqnarray}{RL?s}
\mathbf{I}_{d} \defeq
\begin{bmatrix}
I_{d,a} &
I_{d,n}
\end{bmatrix}^{\text{T}}.
\end{IEEEeqnarray}
Conservation of current through the load is,
\begin{IEEEeqnarray}{RL?s}
 \mathbf{1}^\text{T} \mathbf{I}_{d} = 0. \label{eq_kcl_load}
\end{IEEEeqnarray}
The power flowing through each of the wires to the load is,
\begin{IEEEeqnarray}{RL?s}
\mathbf{S}_{d} 
\defeq
\begin{bmatrix}
S_{d,aa} &
S_{d,nn}
\end{bmatrix}^{\text{T}}
\defeq
\mathbf{U}_{i}
\circ
\mathbf{I}_{d}^* \label{eq_power_def_load_bus}.
\end{IEEEeqnarray}
To force the load value (i.e. constant-power set point), we account for the voltage drop over the load, i.e. $U_{i,a} - U_{i,n}$, 
\begin{IEEEeqnarray}{C}
{S}_{d}^{\text{ref}} = 
(U_{i,a} - U_{i,n}) {I}^*_{d,a} . \label{eq_load_setpoint_def}
\end{IEEEeqnarray}
We now continue to expand the RHS, and use \eqref{eq_kcl_load} and \eqref{eq_power_def_load_bus},
\begin{IEEEeqnarray}{C}
U_{i,a} {I}^*_{d,a} - U_{i,n} {I}^*_{d,a} = U_{i,a} {I}^*_{d,a} + U_{i,n} {I}^*_{d,n} \nonumber\\
= S_{d,aa} + S_{d,nn}  = \mathbf{1}^\text{T}\mathbf{S}_{d}
\end{IEEEeqnarray}
and finally obtain,
\begin{IEEEeqnarray}{C}
{S}_{d}^{\text{ref}} = \mathbf{1}^\text{T}\mathbf{S}_{d} .\label{eq_load_set_point}
\end{IEEEeqnarray}
Note that this does not imply conservation of current,
\begin{IEEEeqnarray}{C}
\eqref{eq_load_set_point} \centernot\implies {I}_{d,a} + {I}_{d,n} = 0.
\end{IEEEeqnarray}
The effective dispatch value for a generator is, 
\begin{IEEEeqnarray}{C}
{S}_{g}^{\text{disp}} \defeq {P}_{g}^{\text{disp}} + j {Q}_{g}^{\text{disp}} = \mathbf{1}^\text{T}\mathbf{S}_{g}. \label{eq_gen_set_point}
\end{IEEEeqnarray}
Finally, KCL at bus $i$ is,
\begin{IEEEeqnarray}{C}
\sum_{lij} \mathbf{I}_{lij} +\sum_{di} \mathbf{I}_{d} - \sum_{gi} \mathbf{I}_{g} = 0 .\label{eq_bus_kcl_iv}
\end{IEEEeqnarray}

\subsection{Canonical Power-Voltage (\svrone)}
We state Ohm's law \eqref{eq_ohms_impedance} in admittance form,
\begin{IEEEeqnarray}{RL?s}
\mathbf{I}_{lij} = \mathbf{Y}_l(\mathbf{U}_i-\mathbf{U}_j)
\end{IEEEeqnarray}  
and substitute this into \eqref{eq_line_power_def} to obtain power flow vector into the branch $l$ in the direction of $i$ to $j$,
\begin{IEEEeqnarray}{RL?s}
\mathbf{S}_{lij} = \mathbf{U}_i \circ\left(\mathbf{Y}_l (\mathbf{U}_i-\mathbf{U}_j)  \right)^* \label{eq_ohms_sv}.
\end{IEEEeqnarray}  
We lift KCL \eqref{eq_bus_kcl_iv} to power variables (vectors with the same sizes as the original current vectors in \eqref{eq_bus_kcl_iv}), by taking the conjugate and multiplying element-wise with $\mathbf{U}_i$,
\begin{IEEEeqnarray}{C}
 \sum_{lij} \mathbf{S}_{lij} +\sum_{di} \mathbf{S}_{d} - \sum_{gi} \mathbf{S}_{g} = 0 \label{eq_bus_kcl_sv}.
\end{IEEEeqnarray}
When we allow  0\,V entries in $\mathbf{U}_i$, the lifted KCL is a relaxation, as it does not imply current-based KCL \eqref{eq_bus_kcl_iv}, due to (possible) multiplication with 0\,V.
This is equivalent to adding a slack variable  to the original KCL, e.g. when $U_{i,n} =$ 0\,V:
\begin{IEEEeqnarray}{C}
\sum_{lij} {I}_{lij,n} +\sum_{di} {I}_{d,n} - \sum_{gi} {I}_{g,n} + I^{\text{slack}}_{i,n} = 0 .\label{eq_bus_kcl_iv_slack}
\end{IEEEeqnarray}
As we are working with complex values, $I^{\text{slack}}_i$ actually represents two degrees of freedom, which can be exploited by optimization solvers. 

\subsection{Power-Voltage Extension (\svrtwo) }
Even though \svrone{} has proven to work well to optimize three-phase Kron-reduced four-wire networks with \emph{only} wye loads\footnote{In this case conservation of current can be arbitrarily satisfied, as the resulting neutral current residual can flow into the ground unimpeded. }, we need to strengthen the feasible set when we model loads/generators between phases (`delta') or phase-to-neutral\footnote{Neutral that is not exactly 0 V.}.
Constraints \eqref{eq_load_set_point},\eqref{eq_gen_set_point} are insufficient in general, as they do not imply  \eqref{eq_kcl_gen} or  \eqref{eq_kcl_load}.
To resolve this, we use the variables $\mathbf{I}_{d},\mathbf{I}_{g}$ and add the constraints  \eqref{eq_kcl_gen}, \eqref{eq_gen_power_def}, \eqref{eq_kcl_load},  \eqref{eq_power_def_load_bus}.

\subsection{Canonical Power-Lifted-Voltage BFM (\swrone)}
We define lifted variables for voltage, $\mathbf{W}_i$, current, $\mathbf{L}_l$, and power, $\mathbf{\bar{S}}_{lij}$,
\begin{IEEEeqnarray}{RL?s}
\mathbf{W}_i \defeq \mathbf{U}_i (\mathbf{U}_i)^{\text{H}}, \label{eq_u_square_def} 
\mathbf{L}_l \defeq \mathbf{I}_{lij} (\mathbf{I}_{lij})^{\text{H}} , \label{eq_i_square_def} 
\mathbf{\bar{S}}_{lij} \defeq \mathbf{U}_i (\mathbf{I}_{lij})^{\text{H}}, \label{eq_s_square_def}
\end{IEEEeqnarray} 
where superscript `H' indicates the conjugate transpose. 
Note that the substitutions don't destroy any information\footnote{Up to a complex scale factor, see \cite{8973484}.}, as long as $\mathbf{U}_i, \mathbf{I}_{lij}$ are not the zero vector - individual entries being zero is fine. 
The power variable $\mathbf{\bar{S}}_{lij}$ is a square matrix, and therefore distinct from the vector $\mathbf{S}_{lij}$, however they are related,
\begin{IEEEeqnarray}{RL?s}
\diag(\mathbf{\bar{S}}_{lij})  = \mathbf{S}_{lij} . \label{eq_power_diag_link}
\end{IEEEeqnarray} 
We can construct a block matrix based on the lifted variables,
\begin{IEEEeqnarray}{RL?s}
\mathbf{M}_{lij} \defeq  
\begin{bmatrix}
\mathbf{U}_i \\
\mathbf{I}_{lij}
\end{bmatrix}
\begin{bmatrix}
\mathbf{U}_i \\
\mathbf{I}_{lij}
\end{bmatrix}^{\text{H}}
=
\begin{bmatrix}
\mathbf{W}_i & \mathbf{\bar{S}}_{lij}\\
(\mathbf{\bar{S}}_{lij})^{\text{H}}& \mathbf{L}_l 
\end{bmatrix}.
\end{IEEEeqnarray} 
In convex relaxation schemes, this constraint gets cast as, 
\begin{IEEEeqnarray}{RL?s}
\rank(\mathbf{M}_{lij} ) = 1, \mathbf{M}_{lij} \succeq 0,
\end{IEEEeqnarray} 
and the rank constraint is dropped to obtain a semi-definite programming formulation.
Ohm's law \eqref{eq_ohms_impedance} is multiplied with its own conjugate transpose to obtain,
\begin{IEEEeqnarray}{RL?s}
\mathbf{W}_j = \mathbf{W}_i - \mathbf{\bar{S}}_{lij}(\mathbf{Z}_l )^{\text{H}} - \mathbf{Z}_l(\mathbf{\bar{S}}_{lij})^{\text{H}} + \mathbf{Z}_l \mathbf{L}_l (\mathbf{Z}_l )^{\text{H}} . \label{eq_ohms_sw}
\end{IEEEeqnarray} 
Furthermore, by multiplying Ohm's law \eqref{eq_ohms_impedance} on the right with $\mathbf{I}_{lij}^{\text{H}}$ and using \eqref{eq_kcl_branch}, we obtain,
\begin{IEEEeqnarray}{RL?s}
\mathbf{\bar{S}}_{lij} + \mathbf{\bar{S}}_{lji} = \mathbf{Z}_l \mathbf{L}_l. \label{eq_flow_loss_sw}
\end{IEEEeqnarray} 
Combining \eqref{eq_bus_kcl_sv} with \eqref{eq_power_diag_link}, we have the BFM formulation with a diagonal KCL expression, as proposed by Gan and Low \cite{gan2014}.

\subsection{Power-Lifted-Voltage BFM Extension (\swrtwo)}
We now generalize the definitions of the load and generator power to matrix variables as well\footnote{Matrix generalizations add a bar above the original symbol.},
\begin{IEEEeqnarray}{RL?s}
\mathbf{\bar{S}}_{d} 
\defeq
\mathbf{U}_{i}
\mathbf{I}_{d}^{\text{H}}
\defeq
\begin{bmatrix}
S_{d,aa} & S_{d,an} \\
S_{d,na}  & S_{d,nn}
\end{bmatrix}
=
\begin{bmatrix}
U_{i,a} I^*_{d,a}  & U_{i,a} I^*_{d,n} \\
U_{i,n} I^*_{d,a}  & U_{i,n} I^*_{d,n}
\end{bmatrix}.
\label{eq_power_def_load_bus_sw}
\end{IEEEeqnarray}
We can partition based on the rows,
\begin{IEEEeqnarray}{RL?s}
\mathbf{\bar{S}}_{d} 
\defeq 
\begin{bmatrix}
\mathbf{S}_{d,a\cdot}  \\
\mathbf{S}_{d,n\cdot}
\end{bmatrix},
\end{IEEEeqnarray}
and recognize that \eqref{eq_kcl_load} implies, 
\begin{IEEEeqnarray}{RL?s}
 \mathbf{1}^\text{T} \mathbf{S}_{d,a\cdot}  = 0, \quad 
 \mathbf{1}^\text{T} \mathbf{S}_{d,n\cdot} = 0. \label{eq_sw_kcl_load}
\end{IEEEeqnarray}
Similarly, for generators we obtain,
\begin{IEEEeqnarray}{RL?s}
\mathbf{\bar{S}}_{g} 
\defeq
\mathbf{U}_{i}
\mathbf{I}_{g}^{\text{H}}
\defeq
\begin{bmatrix}
\mathbf{S}_{g,a\cdot}  \\
\mathbf{S}_{g,n\cdot}
\end{bmatrix}, \,\,\,\,
 \mathbf{1}^\text{T} \mathbf{S}_{g,a\cdot}  = 0, \,\,\,\,
 \mathbf{1}^\text{T} \mathbf{S}_{g,n\cdot} = 0.
\label{eq_power_def_gen_bus_sw}
\end{IEEEeqnarray}
We note valid equalities of the matrix variables w.r.t. the vector-style ones previously defined,
\begin{IEEEeqnarray}{RL?s}
\diag(\mathbf{\bar{S}}_{d})  = \mathbf{S}_{d}, \quad \diag(\mathbf{\bar{S}}_{g})  = \mathbf{S}_{g}.
\end{IEEEeqnarray}
Note that we do \emph{not} project KCL on the rows of $\mathbf{\bar{S}}_{lij}$ or $ \mathbf{L}_{l}$ in a similar fashion to \eqref{eq_sw_kcl_load}-\eqref{eq_power_def_gen_bus_sw}. 
Due capacitive coupling to earth as well as grounding of the neutral, the current through a branch does not in general sum to zero. 
Finally we generalize KCL, using the matrix variables $\mathbf{\bar{S}}$, 
\begin{IEEEeqnarray}{C}
\sum_{lij} \mathbf{\bar{S}}_{lij} +\sum_{di} \mathbf{\bar{S}}_{d} - \sum_{gi} \mathbf{\bar{S}}_{g} = 0 \label{eq_bus_kcl_sw}.
\end{IEEEeqnarray}









\subsection{Voltage Bounds}
For visualizing the feasible sets we enforce voltage bounds,
\begin{IEEEeqnarray}{RL?s}
\underbrace{
\begin{bmatrix}
(U^{\text{min}}_{i,a})^2\\
(U^{\text{min}}_{i,n})^2
\end{bmatrix}
}_{\mathbf{U}^{\text{min}}_i \circ \mathbf{U}^{\text{min}}_i}
\leq
\underbrace{
\begin{bmatrix}
|U_{i,a}|^2\\
|U_{i,n}|^2
\end{bmatrix}
}_{\mathbf{U}_i \circ \mathbf{U}^*_i= \diag(\mathbf{W}_i)}
\leq
\underbrace{
\begin{bmatrix}
(U^{\text{max}}_{i,a})^2\\
(U^{\text{max}}_{i,n})^2
\end{bmatrix}
}_{\mathbf{U}^{\text{max}}_i \circ \mathbf{U}^{\text{max}}_i} \label{eq_voltage_bounds}
\end{IEEEeqnarray}
By setting $U^{\text{min}}_{i,a}$ for instance to 0.9 pu, we typically exclude the low-voltage solution to quadratic equality \eqref{eq_load_setpoint_def}.

\subsection{Feasible sets and generalizations}
Table \ref{tab:feasiblesets} summarizes the variables and constraints for five optimization formulation variants, using the previously-defined constraints.
\begin{table*}[tb]
\centering
\caption{Feasible sets for unbalanced OPF formulation variants. In the objective, $\theta$ is a parameter.} \label{tab:feasiblesets}
\begin{tabular}{@{\extracolsep{4pt}} l  c c c c c} 
\hline
Formulation &  \ivr &  \svrone &  \svrtwo &  \swrone  &  \swrtwo\\ \cline{2-2} \cline{3-4}  \cline{5-6}
& BFM & \multicolumn{2}{c}{BIM } & \multicolumn{2}{c}{BFM } \\  \cline{2-4} 
Math. structure &  \multicolumn{3}{c}{nonconvex QCQP} &  \multicolumn{2}{c}{SDP  }\\ 
Solver &  \multicolumn{3}{c}{\textsc{Ipopt} \cite{ipopt}} &  \multicolumn{2}{c}{\textsc{Clarabel} \cite{goulart_chen} } \\ 
\hline
Objective & \multicolumn{5}{c}{$\min \cos(\theta)\Re({S}_{g}^{\text{disp}} ) + \sin(\theta)\Im({S}_{g}^{\text{disp}} )$ }\\ 
Variables & $\mathbf{U}_i, \mathbf{I}_{lij}, \mathbf{I}_{d}, \mathbf{I}_{g} $  & $\mathbf{U}_i, \mathbf{S}_{lij}, \mathbf{S}_{g} $ & $\mathbf{U}_i, \mathbf{S}_{lij}, \mathbf{S}_{g}$ & $\mathbf{W}_i, \mathbf{L}_l, \mathbf{\bar{S}}_{lij} $ & $\mathbf{W}_i, \mathbf{L}_l, \mathbf{\bar{S}}_{lij} $\\ 
  & & & $\mathbf{I}_{d}, \mathbf{I}_{g}$ & $\mathbf{{S}}_{d} , \mathbf{{S}}_{g} $&  $\mathbf{\bar{S}}_{d}, \mathbf{\bar{S}}_{g} $\\ 
Bus KCL & \eqref{eq_bus_kcl_iv}  & \eqref{eq_bus_kcl_sv} &\eqref{eq_bus_kcl_sv} & \eqref{eq_bus_kcl_sv} &\eqref{eq_bus_kcl_sw}   \\
Ohm's law & \eqref{eq_ohms_impedance}  &  \eqref{eq_ohms_sv} & \eqref{eq_ohms_sv}& \eqref{eq_ohms_sw} & \eqref{eq_ohms_sw} \\
Branch power  & -& -& - & $\mathbf{M}_{lij} \succeq 0$, \eqref{eq_flow_loss_sw} & $\mathbf{M}_{lij} \succeq 0$, \eqref{eq_flow_loss_sw} \\
Load power &\eqref{eq_power_def_load_bus} & - & \eqref{eq_power_def_load_bus} & - & - \\ 
Load set point & \eqref{eq_load_set_point} & \eqref{eq_load_set_point} & \eqref{eq_load_set_point} & \eqref{eq_load_set_point} & \eqref{eq_load_set_point}\\
Load current & \eqref{eq_kcl_load} & - & \eqref{eq_kcl_load} & - & \eqref{eq_sw_kcl_load} \\
Generator power & \eqref{eq_gen_power_def} & - & \eqref{eq_gen_power_def} & - & -\\
Generator dispatch & \eqref{eq_gen_set_point} & \eqref{eq_gen_set_point} & \eqref{eq_gen_set_point} & \eqref{eq_gen_set_point} & \eqref{eq_gen_set_point}\\
Generator current & \eqref{eq_kcl_gen} & - & \eqref{eq_kcl_gen} & - & \eqref{eq_sw_kcl_load} \\
Voltage bounds & \eqref{eq_voltage_bounds} & \eqref{eq_voltage_bounds}  & \eqref{eq_voltage_bounds} & \eqref{eq_voltage_bounds}  & \eqref{eq_voltage_bounds} \\
 \hline
\end{tabular}
\end{table*}
Implementation of these equations in nonlinear optimization frameworks typically requires translating the complex-valued expressions to sets of equivalent real-valued ones.
Furthermore, lines should be modelled with shunt admittances, typically used to represent capacitive behavior w.r.t earth. 
Note that in general grid codes specify phase-to-neutral voltage magnitude limits, not phase-to-ground. 
Geth and Ergun  \cite{geth2021a}  present the corresponding real-value expressions, and derive expressions for a variety of bounds.
Claeys et al. \cite{claeys2020a} propose a framework to represent unbalanced transformers across variable spaces.

\section{OPF for a 2-bus 2-wire case}
Note that extensive numerical results showcasing 0\,V issues ($\geq 128$ networks, up to 2156 nodes) have been published in papers \cite{vanacker2022,claeys2022} co-authored  by this paper's author. 
Therefore, in this work, the focus is to provide \emph{insight}, using  a 2-bus test case (data in Table \ref{tab:2bustestcase}). 
The problems are implemented in \textsc{JuMP} \cite{DunningHuchetteLubin2017}. 
The nonconvex formulations are solved with \textsc{Ipopt} \cite{ipopt}; the semidefinite relaxation is solved with \textsc{Clarabel} \cite{goulart_chen}. 
We solve repeatedly for different values  $0\leq \theta \leq 2\pi$ in the objective, to sample the edge of the solution space. 

\begin{table}[tb]
\centering
\caption{Data for 2-bus test case (per unit)} \label{tab:2bustestcase}
\begin{tabular}{l  c c c c c} 
\hline
Voltage source setup & $U_{i,a} = 1\angle 0 $, $U_{i,n} = 0 $\\
Line resistance & $\mathbf{R}_l$ = $0.05\cdot [1\;0.1;\;0.1\; 1]$ \\
Line reactance & $\mathbf{X}_l$ = $0.04 \cdot[1\; 0.5;\;0.5\; 1]$ \\
Load set point & ${S}_{d}^{\text{ref}} = 1+j0.5$ \\
Voltage  bounds & $\mathbf{U}^{\text{min}}_i = [0.9 \,\, \textcolor{red}{0}]^{\text{T}}$ ,
$\mathbf{U}^{\text{max}}_i = [1.1\,\, 1.1]^{\text{T}}$ \\
 \hline
\end{tabular}
\end{table}

\begin{figure}[tbh]
  \centering
    \includegraphics[width=0.75\columnwidth]{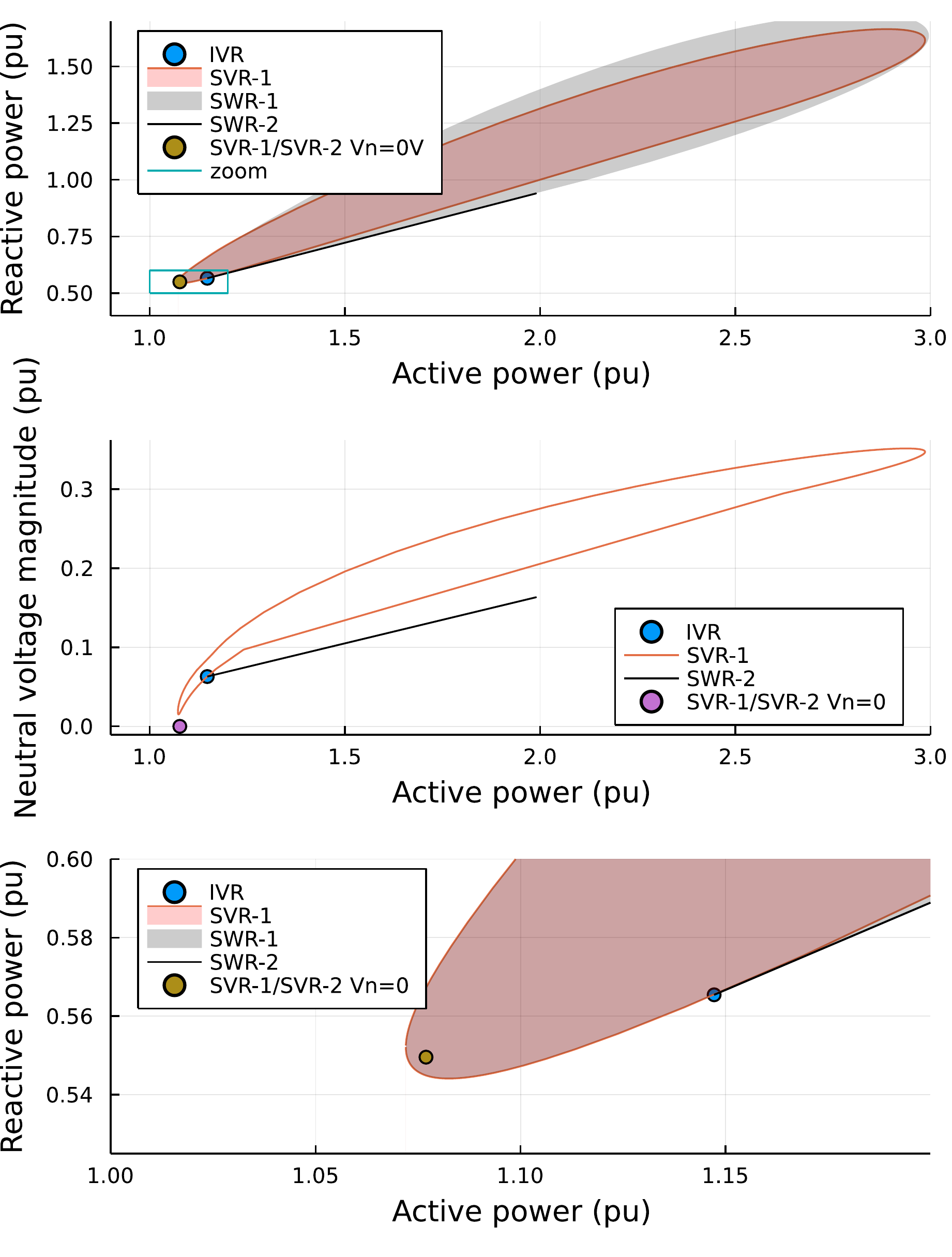}
  \caption{The top figure illustrates the solution sets for the formulations, the bottom one zooms into the solutions with minimum active power. The middle figure indicates the corresponding value of the neutral voltage magnitude at bus $j$ for the solutions on the edges of the solution sets.}  \label{fig_feasib}
\end{figure}

Note that this is actually a power flow simulation, i.e. feasibility, problem, as with no physical degree of freedom, as the generator needs to produce a specific amount of power to satisfy the demand, taking into account the branch losses.
As will be shown though, only  \ivr{} has a single (high-voltage) solution, all other forms have multiple solutions. 
The solutions sets are visualized in Fig.~\ref{fig_feasib} and are discussed next.

\subsection{Numerical results: amount of solutions}

\subsubsection{\ivr{} has a single solution}
The solution to this problem is uniquely specified by the voltage phasor at bus $j$, i.e.,
\begin{IEEEeqnarray}{RL?s}
\mathbf{U}_j = 
\begin{bmatrix}
 0.937066 + j0.002500 \\
 0.062934 - j0.002500
\end{bmatrix}, 
\mathbf{U}_j^{\text{mag}} = 
\begin{bmatrix}
 0.937069 \\
 0.062983
\end{bmatrix}, 
\end{IEEEeqnarray}
and a phase-to-neutral voltage magnitude of $0.874146$ pu.
The generator output is,
$  {S}_{g}^{\text{disp}} = 1.147226 + j0.565434 $.
This solution is that of the circuit in Fig.~\ref{fig_SU-IV}.


\subsubsection{\svrtwo{} has two solutions}
The \ivr{} solution is feasible w.r.t. \svrtwo, however \svrtwo{} has a second spurious solution, 
\begin{IEEEeqnarray}{RL?s}
\mathbf{U}_j = 
\begin{bmatrix}
0.924699 - j0.007591 \\
                0 + j0
\end{bmatrix}, 
\mathbf{U}_j^{\text{mag}} = 
\begin{bmatrix}
 0.924730 \\
 0
\end{bmatrix}, 
\end{IEEEeqnarray}
associated with generator dispatch, $  {S}_{g}^{\text{disp}} =  1.076921 + j0.549558  $.
The relaxation gap of this solution is $6.1282\%$.
The voltage drop across the load differs, with the phase-to-neutral voltage magnitude at $ 0.924$ where it was $0.874$ pu for \ivr{}.
With the neutral being hard grounded at both buses, this solution is equivalent to that of a transformed circuit with Kron's reduction of the neutral ($\mathbf{Z}^{\text{Kr}}_l = 0.052622 + j0.033902$pu). 

\begin{figure}[tbh]
  \centering
    \includegraphics[width=0.65\columnwidth]{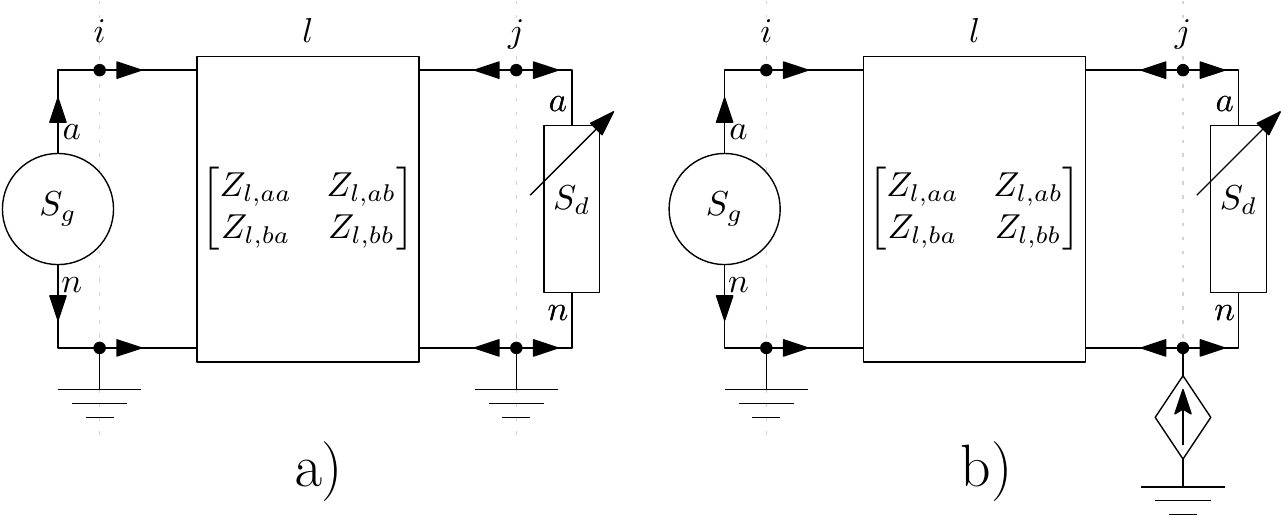}
  \caption{a) relaxed solution allowed by \svrtwo{} formulation; b) relaxed solutions allowed by \svrone{} formulation.}  
  \label{fig_SU_solutions_grounded}
\end{figure}

This spurious solution is that of the circuit in Fig.~\ref{fig_SU_solutions_grounded}a, where the neutral at bus $j$ is grounded.
The existence of two distinct solutions suggests the problem for bigger networks has combinatorial (discontinuous) features, which may limit the performance of many algorithms, including the NLP interior-point algorithm in \textsc{Ipopt}.

\subsubsection{\svrone{} has infinite solutions}
A two-dimensional infinite set of solutions is feasible w.r.t. \svrone{}. 
Note that the \ivr{} solution is on the edge of this set, whereas the spurious \svrtwo{} solution is on the interior. 
Despite the \ivr{} solution being on the edge, it is not in an ideal location, as it does not align with the minimum active power search direction, so likely leads to nonzero relaxation gap under a broad range of objectives.

For the 2-bus network, the combination of dropping the load KCL and allowing the 0 neutral voltage solution is equivalent to the solutions of the circuit presented in Fig.~\ref{fig_SU_solutions_grounded}b, where a variable complex current source is added to the bus $j$.
The $\min {P}_{g}^{\text{disp}}$ solution is,
\begin{IEEEeqnarray}{RL?s}
\mathbf{U}_j = 
\begin{bmatrix}
    \phantom{-}0.926241 - j 0.013080 \\
 -0.007515 - j 0.016886
\end{bmatrix}, 
\mathbf{U}_j^{\text{mag}} = 
\begin{bmatrix}
 0.926333 \\
 0.018483
\end{bmatrix}, 
\end{IEEEeqnarray}
with phase-to-neutral voltage magnitude at $0.933764$ pu.
The generator output is $  {S}_{g}^{\text{disp}} =  1.071996 + j 0.552133  $, i.e. a gap of $6.5575 \%$

\subsubsection{\swrtwo{} has infinite solutions}
The solution space is a 1-dimensional infinite set.
The \ivr{} solution is an extreme point on this line, in the direction of minimum apparent power, therefore  leading to a zero relaxation gap across a broad range of real-world objectives.  
The 0\,V \svrtwo{} solution is \emph{not} in the feasible set of \swrtwo{}.

\subsubsection{\swrone{} has infinite solutions}
The solution space is a 2-dimensional infinite set, with a nonzero gap. 
We do not establish if it is comparable with \svrone{}. 
The $\min {P}_{g}^{\text{disp}}$ solution is the same as that of \svrone{}.

\subsection{Discussion}
The exclusion of the 0\,V solution in \swrtwo{} is due to the combination of 1) KCL as a matrix equality 2) conservation of current through loads/generators. Adding either constraint individually to \swrone{} is not sufficient to obtain a 0 gap.
That suggests the spurious solution in \svrtwo{} can be excluded by also using the matrix KCL \eqref{eq_bus_kcl_sw} and generalizing \eqref{eq_ohms_sv} to the outer product (instead of element-wise).
Alternatively, as voltages are inherently relative, we can also choose to redefine the ground voltage value, at for instance 1\,V.

Modelling loads with a nonzero voltage at both terminals requires additional constraints. 
In \svrtwo{}, the current variables could still be eliminated in favor of power variables (i.e. as $(S/U)^*$), following the approach in \cite{claeys2022}, but this may result in higher order polynomial expressions. 
Furthermore, constant-power three-phase loads  have multiple solutions, with the amount depending on the configuration and neutral grounding impedance~\cite{araujo2016,Wang03}. 
Deriving valid bounds to  exclude the `low-voltage solution' may speed up the computation.


\section{Conclusions}
In the absence of voltage magnitude lower bounds, the `exact' canonical unbalanced power-voltage formulation  is a relaxation of the circuit laws with a nonzero gap. 
Conversely, the unbalanced SDP BFM formulation - often trivialized for its theoretical inexactness,  has a zero gap under the same conditions. 
The nonzero gap solution permitted by the power-voltage formulation, itself characterized by a zero voltage on one of the bus terminals, can furthermore be proven infeasible by the SDP relaxation.
In what amounts to a reversal of fortune w.r.t. transmission-style (positive-sequence) OPF in terms of reliability, convergence and gap, with SDP now performing better than NLP, bound tightening based on a convex relaxation could be performed to  generate valid lower bounds for the power-voltage form, thereby making the solution unique again.

\bibliographystyle{IEEEtran}

\begin{thebibliography}{10}
\providecommand{\url}[1]{#1}
\csname url@samestyle\endcsname
\providecommand{\newblock}{\relax}
\providecommand{\bibinfo}[2]{#2}
\providecommand{\BIBentrySTDinterwordspacing}{\spaceskip=0pt\relax}
\providecommand{\BIBentryALTinterwordstretchfactor}{4}
\providecommand{\BIBentryALTinterwordspacing}{\spaceskip=\fontdimen2\font plus
\BIBentryALTinterwordstretchfactor\fontdimen3\font minus
  \fontdimen4\font\relax}
\providecommand{\BIBforeignlanguage}[2]{{%
\expandafter\ifx\csname l@#1\endcsname\relax
\typeout{** WARNING: IEEEtran.bst: No hyphenation pattern has been}%
\typeout{** loaded for the language `#1'. Using the pattern for}%
\typeout{** the default language instead.}%
\else
\language=\csname l@#1\endcsname
\fi
#2}}
\providecommand{\BIBdecl}{\relax}
\BIBdecl

\bibitem{kardos2018complete}
J.~Kardoš, D.~Kourounis, O.~Schenk, and R.~Zimmerman, ``Complete results for a
  numerical evaluation of interior point solvers for large-scale optimal power
  flow problems,'' \emph{arXiv:1807.03964}, 2018.

\bibitem{ipopt}
A.~{W\"{a}chter} and L.~{Biegler}, ``On the implementation of an interior-point
  filter line-search algorithm for large-scale nonlinear programming,''
  \emph{Math. Program.}, vol. 106, no.~1, pp. 25--57, 2006.

\bibitem{5491276}
R.~D. Zimmerman, C.~E. Murillo-Sánchez, and R.~J. Thomas, ``Matpower:
  Steady-state operations, planning, and analysis tools for power systems
  research and education,'' \emph{IEEE Trans. Power Syst.}, vol.~26, no.~1, pp.
  12--19, 2011.

\bibitem{KARDOS2022108613}
J.~Kardoš, D.~Kourounis, O.~Schenk, and R.~Zimmerman, ``{BELTISTOS}: A robust
  interior point method for large-scale optimal power flow problems,''
  \emph{Elect. Power Syst. Res.}, vol. 212, p. 108613, 2022.

\bibitem{vanacker2015}
T.~Van~Acker, D.~Van~Hertem, D.~Bekaert, K.~Karoui, and C.~Merckx,
  ``Implementation of bus bar switching and short circuit constraints in
  optimal power flow problems,'' in \emph{IEEE PowerTech}, 2015.

\bibitem{en14082160}
A.~K. Barnes, J.~E. Tabarez, A.~Mate, and R.~W. Bent, ``Optimization-based
  formulations for short-circuit studies with inverter-interfaced generation in
  {PowerModelsProtection.jl},'' \emph{Energies}, vol.~14, no.~8, 2021.

\bibitem{7587824}
T.~Mühlpfordt, T.~Faulwasser, and V.~Hagenmeyer, ``Solving stochastic ac power
  flow via polynomial chaos expansion,'' in \emph{IEEE Conf. Control App.},
  2016, pp. 70--76.

\bibitem{vanacker2022}
T.~Van~Acker, F.~Geth, A.~Koirala, and H.~Ergun, ``General polynomial chaos in
  the current–voltage formulation of the optimal power flow problem,''
  \emph{Elect. Power Syst. Res.}, vol. 211, p. 108472, 10 2022.

\bibitem{geth2022}
F.~Geth and T.~Van~Acker, ``Harmonic optimal power flow with transformer
  excitation,'' \emph{Elect. Power Syst. Res.}, vol. 213, p. 108604, 12 2022.

\bibitem{ARAUJO2013632}
L.~R. de~Araujo, D.~R.~R. Penido, and F.~de~Alcântara~Vieira, ``A multiphase
  optimal power flow algorithm for unbalanced distribution systems,''
  \emph{Int. J. Elect. Power Energy Syst.}, vol.~53, pp. 632--642, 2013.

\bibitem{claeys2022}
S.~Claeys, F.~Geth, and G.~Deconinck, ``Optimal power flow in four-wire
  distribution networks: Formulation and benchmarking,'' \emph{Elect. Power
  Syst. Res.}, vol. 213, p. 108522, 12 2022.

\bibitem{chandra2022}
\BIBentryALTinterwordspacing
C.~K. Jat, J.~Dave, D.~Van~Hertem, and H.~Ergun, ``Unbalanced opf modelling for
  mixed monopolar and bipolar hvdc grid configurations,'' 2022. [Online].
  Available: \url{https://arxiv.org/abs/2211.06283}
\BIBentrySTDinterwordspacing

\bibitem{gan2014}
L.~Gan and S.~H. Low, ``Convex relaxations and linear approximation for optimal
  power flow in multiphase radial networks,'' in \emph{Proc. Power Syst. Comp.
  Conf.}, 2014.

\bibitem{geth2021a}
F.~Geth and H.~Ergun, ``Real-value power-voltage formulations of, and bounds
  for, three-wire unbalanced optimal power flow,'' \emph{arXiv:2106.06186}, pp.
  1--12, 2021.

\bibitem{8973484}
F.~Geth and C.~Coffrin, ``Direct method to recover current and voltage in
  multi-conductor optimal power flow models,'' in \emph{IEEE Power Energy Soc.
  General Meeting}, 2019, pp. 1--5.

\bibitem{goulart_chen}
\BIBentryALTinterwordspacing
P.~Goulart and Y.~Chen. [Online]. Available:
  \url{https://oxfordcontrol.github.io/ClarabelDocs/stable/}
\BIBentrySTDinterwordspacing

\bibitem{claeys2020a}
S.~Claeys, G.~Deconinck, and F.~Geth, ``Decomposition of n-winding transformers
  for unbalanced optimal power flow,'' \emph{IET Gen. Trans. Distrib.},
  vol.~14, pp. 5816--5822, 12 2020.

\bibitem{DunningHuchetteLubin2017}
I.~Dunning, J.~Huchette, and M.~Lubin, ``{JuMP}: A modeling language for
  mathematical optimization,'' \emph{SIAM Rev.}, vol.~59, no.~2, pp. 295--320,
  2017.

\bibitem{araujo2016}
L.~R. de~Araujo, D.~R.~R. Penido, S.~Carneiro, and J.~L.~R. Pereira, ``A study
  of neutral conductors and grounding impacts on the load-flow solutions of
  unbalanced distribution systems,'' \emph{IEEE Trans. Power Syst.}, vol.~31,
  pp. 3684--3692, 9 2016.

\bibitem{Wang03}
Y.~Wang and W.~Xu, ``The existence of multiple power flow solutions in
  unbalanced three-phase circuits,'' \emph{IEEE Trans. Power Syst.}, vol.~18,
  no.~2, pp. 605--610, 2003.

\end{thebibliography}

\end{document}